\documentclass[12pt]{JHEP3}
\usepackage{mathrsfs}
\usepackage{amsmath}
\usepackage{epsfig}
\title{ The fractal dimension of the Riemann zeta zeros}
\author{Jingbo Wang \\
  Department of Applied Physics,
 Xi'an Jiaotong University,
 Xi'an, 710049, People's Republic of China\\
 \email{ shuijing31@gmail.com }}

\abstract { In this paper, we consider the nontrivial zeros of the
Riemann zeta function as the eigenvalues of the Dirac operator on a
fractal manifold. From the heat kernel expansion, we figure out that
the fractal dimension of the manifold is about 1.1-1.2. Also we
compare this result to the random matrix theory and the quantum
chaos theory.}

 \preprint{} \dedicated{} \maketitle
\begin{document}
\section{Introduction}
The Riemann hypothesis(RH)\cite{R1} is that all the nontrivial zeros
of the zeta function have real part 1/2, so that the quantities
${t_j}$ defined by
\begin{equation}\label{1}
    \zeta (\frac{1} {2} + it_j ) = 0
\end{equation}
are all real. There are strong evidence to support this hypothesis:
the first billions of zeros are on the line, and also more than
$40\%$ of the nontrivial zeros satisfy the Riemann hypothesis. There
are also mathematical evidence, such as the Deligne's proof of the
corresponding Riemann hypothesis for Zeta function of arbitrary
varieties over finite fields, or the Weil's explicit formula. For a
review, see the paper written for the Millennium by
E.Bombieri\cite{Bo}, or the paper\cite{Co} by J.B.Conrey.

Polya and Hilbert made the conjecture that the imaginary part of the
Riemann zeros could be the eigenvalues of a Hermitian operator. Then
the RH would follow, since Hermitian operator have real eigenvalues.
In the 50's Selberg found a remarkable duality between the length of
the geodesics on a Riemann surface and the eigenvalues of the
Laplacian operator defined on it\cite{Se}. Quite independently of
Selberg's work, Montgomery showed that the Riemann zeros are
distributed randomly and obeying locally the statistical law of the
Random Matrix Theory(RMT)\cite{Mo}. The matrix related to the
Riemann zeros is the gaussian unitary ensemble(GUE) associated to
the random hermitian matrices. Montgomery's results found an
impressive numerical confirmation in the work of Odlyzko in the
80's, so the GUE law is nowadays called the Mongomery-Odlyzko
law\cite{Od}.

A further step along this physical approach to RH was taken by
Berry, who noticed a similarity between the formula yielding the
fluctuations of the number of the zeros around its average position,
and a formula giving the fluctuations of the energy levels of a
Hamiltonian obtained by the quantization of a classical chaotic
system\cite{Be}. The quantum chaos approach can explain the
deviations from the GUE law found numerically by Odlyzko. By this
model suffer an overall sign problem and this problem lead Connes to
propose an abstract approach to the RH based on discrete
mathematical objects known as adeles\cite{Con}. Those approaches
give two possible physical realizations of the Riemann zeros, either
as point like spectra or as missing spectra in a continuum. The next
step came in 1999 when Berry and Keating\cite{BK} on one hand and
Connes\cite{Con} on the other, proposed that the classical
Hamiltonian $H=xp$, where $x$ and $p$ are the position and momenta
of a 1D particle, is closely related to the Riemann zeros. They
choose different regularizations of this Hamiltonion and gave
different realization of the Riemann zeros. And all these
semiclassical results are heuristic and lack so far of a consistent
quantum version. And there are further work along this
line\cite{Si}.

In this paper, we will go along with other approach, that is to
consider the zeros as the eigenvalues of the Dirac operator on a
suitable manifold. With the heat kernel expansion we can get the
dimension of this manifold and find that the manifold is fractal and
has fractal dimension about $1.1-1.2$, closed to the 1 dimension.
And our result don't dependent on the models.
\section{The Riemann Zeta Zeros}
The Riemann Zeta function is defined in the half-plane $Re(s)>1$ by
\begin{equation}\label{2}
    \zeta (s) = \sum\limits_{n = 1}^\infty  {\frac{1} {{n^s }}}  =
\prod\limits_{p(prime)} {\frac{1} {{1 - p^{ - s} }}}.
\end{equation}
and has a meromorphic continuation to $\mathbb{C}$. This function
satisfy the functional equation
\begin{equation}\label{3}
\xi (s): = \pi ^{ - s/2} \Gamma (s/2)\zeta (s) = \xi (1 - s).
    \end{equation}
It is known that the complex zeros of $\zeta(s)$ lie in the
"critical strip" $0<Re(s)<1$, and Riemann hypothesis states that in
fact all these zeros lie on the "critical line" $Re(s)=1/2$, that is
the $t_j$ on \ref{1} are all real. The asymptotic density of the
zeros is given by
\begin{equation}\label{4}
d(t) = \frac{1} {{2\pi }}\log (\frac{t} {{2\pi }}) + O(\frac{1}
{{t^2 }})
\end{equation}
and therefore that the mean spacing between the zeros decreases
logarithmically with increasing $t$. Define the scaled Riemann zeros
so as to have unit mean spacing, that is
\begin{equation}\label{5}
u_j  = \frac{{t_j }} {{2\pi }}\log \frac{{t_j }} {{2\pi }}.
\end{equation}
We will consider those numbers as the eigenvalues of the Dirac
operator on some conjectured manifold, and compare to the results
from the random matrix theory and quantum chaos theory.
\section{The fractal dimension of the Riemann zeros}
The reason why we consider the zeros as the eigenvalues of the Dirac
operator is that, from the functional equation \ref{3} we know that
if $1/2+it_j$ is the zero of the Riemann zeta function, $1/2-it_j$
is also the zero of the function, so $
 \pm t_j$ both are eigenvalues, thus we choose the Dirac operator
 instead of the Laplacian operator.

 We assume that this Dirac operator has the similar heat kernel
 expansion as the usual one, and this assumption can be verified by
 our results. From the heat kernel expansion we can define the
 fractal dimension as follows\cite{Mo}:
 \begin{equation}\label{6}
    D_s  =  - \frac{{d\ln (Trace(e^{ - D^2 /\Lambda ^2 } ))}}{{d\ln
\Lambda }}.
 \end{equation}
Substitute the scaled Riemann zeros ${u_j}$ \cite{Od2} into \ref{6},
and plot the figure of the fractal dimension varying with the
$\Lambda$. \EPSFIGURE{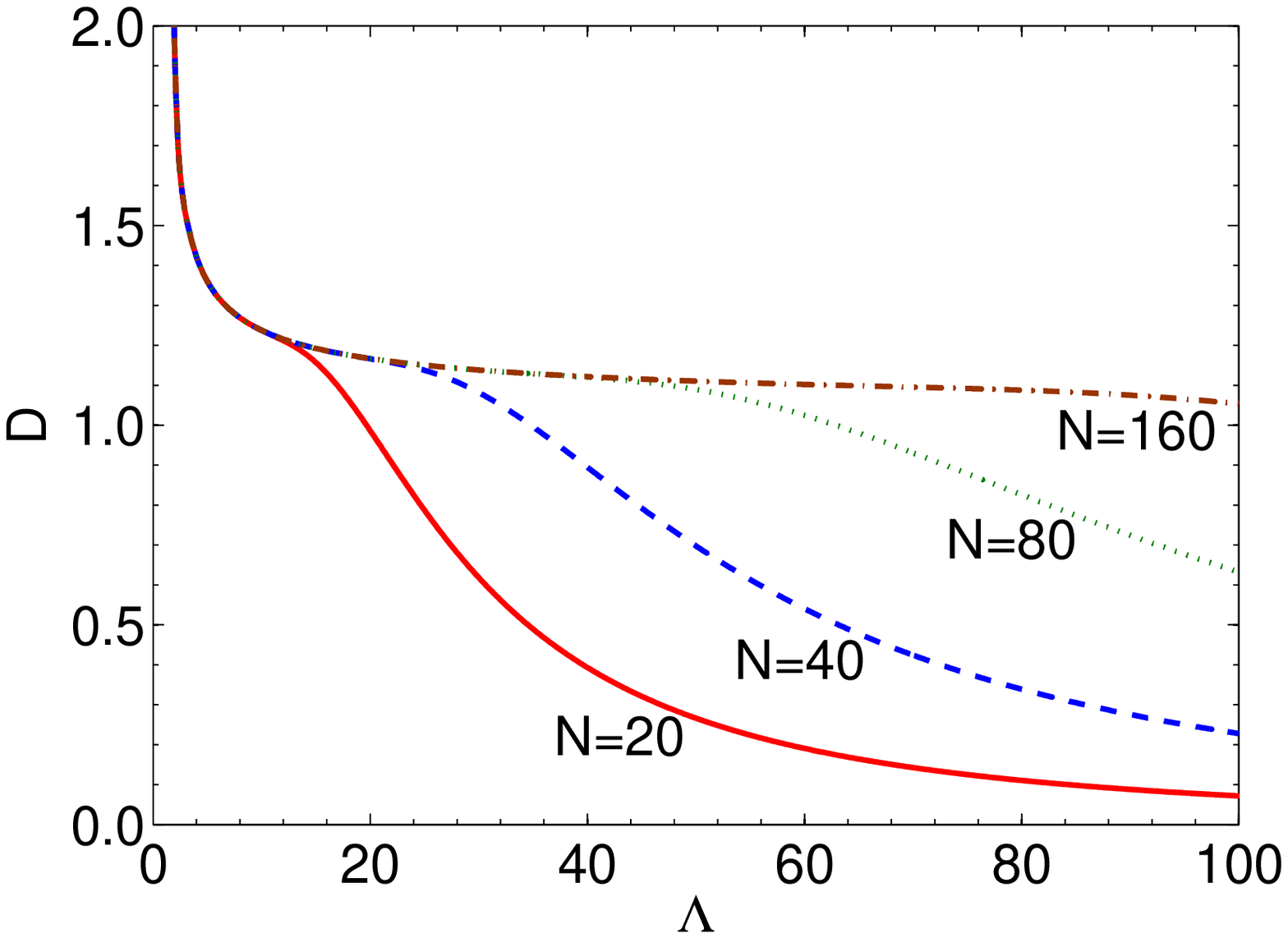} {The fractal dimension of the
Riemann zeta zeros for different numbers} From the above figure we
can see that although the fractal dimension general vary with the
$\Lambda$, at the middle region, the curve is nearly flat, and
correspondent to
\begin{equation}\label{7}
    D_s=1.1-1.2
\end{equation}
In the previous paper\cite{Wang}, we show that for ordinary sphere,
this region give the correct dimension, so here we also consider
this number is the dimension of the underling manifold. We also
notice that for more zeros, the nearly flat region will be longer.

The GUE conjecture state that all statistics for zeta zeros and
eigenvalues of Hermitian matrices match. And the numerical
calculation performed by Odlyzko confirmed this conjecture. But it
can't tell us exactly which Hermitian matrix gives the zeta zeros.
On the other hand, the quantum chaos theory give some clues for the
dynamics behind the zeta function, that is chaotic, time-asymmetry
system. But the system is 1 dimension, and can't give the above
figure. So we think that this approach is unsatisfying.
\section{Conclusion}
We get the fractal dimension from the scaled zeros of the Riemann
zeta function, so it must be obeyed by any system which is used to
model the Riemann zeta function. Unfortunately, the current physical
models, such as the GUE matrices, quantum chaos system don't satisfy
this picture.

The fractal dimension is about $D_s=1.1-1.2$, not away from the
integer 1, but show the fractal character of the Riemann zeta
function. The similarity between the Riemann hypothesis and the
quantum gravity has been point out by Connes and Marcolli in their
book\cite{CM}, and our result confirm this relations.
\acknowledgments{This work was partly done at Beijing Normal
University. This research was supported in part by the Project of
Knowledge Innovation Program (PKIP) of Chinese Academy of Sciences,
Grant No. KJCX2.YW.W10 }

\bibliography{riemann1}
\end{document}